\documentclass[12pt]{article}
\pdfoutput=1
\usepackage{amssymb}
\usepackage{amsmath}
\usepackage{graphicx}
\usepackage{fancybox}
\usepackage{epsfig}
\usepackage{latexsym}
\usepackage[applemac]{inputenc}
\usepackage{ae,aecompl}
\usepackage{amsthm}
\usepackage{mathrsfs}
\usepackage{amsfonts,amssymb}
\usepackage{pstricks}

\renewcommand{\geq}{\geqslant}
\renewcommand{\leq}{\leqslant}

\newcommand{\op}[1]{\operatorname{#1 }}

\newcommand{\R}{\mathbb R}
\newcommand{\N}{\mathbb N}
\newcommand{\Z}{\mathbb Z}

\newtheorem{thm}{Theorem}[section]
\newtheorem{lemma}[thm]{Lemma}

\newtheorem{pro}[thm]{Proposition}
\newtheorem{lem}[thm]{Lemma}
\newtheorem{cor}[thm]{Corollary}
\newtheorem{defi}[thm]{Definition}

\newtheorem{ques}{Question}

\newtheorem{rqe}[thm]{Remark}

\usepackage{enumerate}

\title{On limits of Graphs Sphere Packed in Euclidean Space and Applications}
\date{October 2010}
\author{Itai Benjamini and Nicolas Curien}
\begin{document}
\maketitle

\abstract{The core of this note is the observation that links between circle packings of graphs and potential theory developed in \cite{BeSc01} and \cite{HS} can be extended to higher dimensions. In particular, it is shown that every limit of finite graphs sphere packed in $\R^d$ with a uniformly-chosen root is $d$-parabolic. We then derive few geometric corollaries. E.g.\,every infinite graph packed in $\R^{d}$ has either
strictly positive isoperimetric Cheeger constant or admits arbitrarily large finite sets $W$ with boundary  size which satisfies $ |\partial W| \leq |W|^{\frac{d-1}{d}+o(1)}$. Some open problems and conjectures are gathered at the end.}

\section{Introduction}
The theory of random planar graphs, also known as two-dimensional quantum gravity in the physics literature, has been rapidly developing for the last ten years, see \cite{B10} for a survey. The analogous theory in higher dimension is notoriously hard and not much established so far, this is due in particular to the fact that enumeration techniques and bijective representations are missing, see for instance \cite{BeZi09}. \\
However there are a couple of two dimensional  results that are not depending on enumeration. E.g.\,in \cite{BeSc01}, circle packing theory is used to show that limits (see Definition \ref{limits}) of finite random \emph{planar} graphs of bounded degree with a uniformly-chosen root are almost surely recurrent. The goal of this note is to extend this result into higher dimensions and to draw some consequences and conjectures. \\
We recall that recurrence means that the simple random walk on the graph returns to the origin almost surely, or in a potential theory terminology that the graph is parabolic. A graph is parabolic if and only if it supports no flow with one source of flux $1$, no sinks, and with gradient in  $\mathbb{L}^2$. Replacing $2$ by $d\geq 3$ yields to the concept of $d$-parabolicity, see \cite{S94} and Section \ref{dpot}.\\
The analogous of circle packing theory in dimension $d$ is easy to describe. A graph is sphere packable in $\R^d$ if and only if it is the tangency graph of a collection of $d$-dimensional balls with disjoint interiors: the balls of the packing correspond to the vertices of the graph and the edges to tangent balls, see Section \ref{packings}. The theory of circle packings of planar graphs is well developed and its relation to conformal geometry is well established, see the beautiful survey \cite{R}.  The higher dimensional version is not as neat. First, although all finite planar graphs (without loops nor multiple edges) can be realized as the tangency graph of a circle packing in $\R^2$ (see below), yet there are no natural families of graph packed in $\R^d$ for $d\geq 3$. Second, circle packings relates to $\mathbb{L}^2$ potential theory while in higher dimension the link is to $d$-potential theory, which is less natural and where the probabilistic interpretation is laking. Still useful things can be proved and conjectured.  Indeed the main observation of this note is that links between circle packings of graphs and potential theory over the graph (see \cite{HS}) can be extended to higher dimensions, leading in particular to a generalization of \cite[Theorem 1.1]{BeSc01} and suggests many problems for further research. For a precise formulation of our main theorem (Theorem \ref{dpara}) we must introduce several technical notions and definitions in the coming sections.



We hope that this minor contribution will open the doors for three and higher dimensional theory of sphere packing and quantum gravity. The proofs essentially follow that of \cite{BeSc01} and \cite{HS} with the proper modifications followed by a  report on some new geometric applications. For example we prove under a local bounded geometry assumption defined in the next section that a sequence of $k$-regular graphs with growing girth can not be all  packed in a fixed dimension and that  every infinite graph packed in $\R^{d}$  either has strictly positive isoperimetric Cheeger constant or admits arbitrarily large finite sets $W$ with boundary  size which satisfies $ |\partial W| \leq |W|^{\frac{d-1}{d}+o(1)}$.
\medskip

Note that very recently the isoperimetric criterion of Proposition \ref{geo} was used in \cite{KPP} to prove that acute triangulations of the space $\R^d$ do not exist for $d \geq 5$.

\section{Notations and terminology}
In the following, unless indicated, all graphs are locally finite and connected.
\subsection{Packings}
\label{packings}
\begin{defi} A $d$-dimensional sphere packing or shortly \emph{$d$-sphere packing} is a collection $P = (B_{v}, v \in V )$ of $d$-dimensional
balls of centers $C_v$ and radii $r_{v}>0$ with disjoint interiors in $\R^{d}$. We associated to $P$
an unoriented  graph $G = (V,E)$ called \emph{tangency graph}, where we put an edge between two vertices $u$ and $v$ if and only if the balls $B_{u}$ and $B_{v}$ are tangent.
\end{defi}

An \emph{accumulation point} of a sphere packing $P$ is an accumulation point of the centers of the balls of $P$. Note that the name ``sphere packing'' is unfortunate since it deals with balls. However this terminology is common and we will use it. The $2$-dimensional case is well-understood, thanks to the following Theorem. \begin{thm}[Circle Packing Theorem] A finite graph $G$ is the tangency graph of a $2$-sphere packing if and only if $G$ is planar and contains no multiple edges
nor loops. Moreover if $G$ is a triangulation then this packing is unique up to M\"obius transformations.
\end{thm}
This beautiful result has a long history, we refer to \cite{St05} and \cite{R} for further information. When $d = 3$, very little is known. Although some necessary conditions for a finite
graph to be the tangency graph of a $3$-sphere packing are provided in \cite{KS} (for a related higher dimensional result see \cite{A}), the characterization of $3$-sphere packable graphs is still open (see last section). For packing of infinite graphs
see \cite{BS09}. To bypass the lack of a result similar to the last theorem in dimension $3$ or higher, we will restrict ourselves to packable graphs, that are graphs which admit a sphere packing representation. One useful lemma in circle packing theory is the so-called ``Ring lemma''  that enables us to control the size of tangent circles under a bounded-degree assumption. 
\begin{lemma}[Ring Lemma, \cite{RS}] There is a constant $r>0$ depending only on $n \in \mathbb{Z}_{+}$ such that if $n$ circles surround the unit disk then each circle has radius at least $r$.
\end{lemma} 

Here also, since we have no analogous of the Ring Lemma in high dimensions, we will required an additional property on the packings.

\begin{defi} Let $M > 0$.   A $d$-sphere packing $P=(B_{v},v \in V)$  is \emph{$M$-uniform} if, for any tangent balls $B_{u}$ and $B_{v}$ of radii $r_{u}$ and $r_{v}$ we have
$$ \frac{r_{u}}{r_{v}} \leq M.$$
A graph $G$ is $M$-uniform in dimension $d$, if it is a tangency graph of a $M$-uniform sphere packing in $\R^d$. 
\end{defi}
\begin{rqe} Note that an $M$-uniform graph in dimension $d$ has a maximal degree bounded by a constant depending only on $M$ and $d$.
\end{rqe}
\begin{rqe} By the Ring Lemma, every planar graph without loops nor multiple edges is $M$-uniform in dimension $2$ where $M$ only depends of the maximal degree of the graph. The same holds in dimension $3$ provided that the complex generated by the centers of the spheres is a tetrahedrangulation (that is all simplexes of dimension 3 are tetrahedrons), see \cite{Vasilis}.
\end{rqe}


\subsection{$d$-parabolicity}
\label{dpot}

The classical theory of electrical networks and $2$-potential theory is long studied and well understood, in particular due to the connection with simple random walk (see for example \cite{DoSn84} for a nice introduction). On the other hand, non-linear potential theory is much more complicated and still developing, for background see \cite{S94}. A key concept for $d$-potential theory is the notion of extremal length and its relations with parabolicity (extremal length is common in complex analysis and was imported in the discrete setting by Duffin \cite{Du62}). We present here the basic definitions that we use in the sequel.\\
Let $G=(V,E)$ be a locally finite connected graph. For $v\in V$ we  let $\Gamma(v)$ be the set of all semi-infinite paths in $G$ starting from $v$. If $m : V \to \R_{+}$ is assigning length to vertices, the length of a path $\gamma$ in $G$ :
$$ \op{Length}_{m}(\gamma) := \sum_{v\in \gamma}m(v).$$
If $m \in \mathbb{L}^d(V)$, we denote by $\|m\|_{d}$ the usual $\mathbb{L}^d$ norm $(\sum_v m(v)^d)^{1/d}$. The graph $G$ is \emph{$d$-parabolic} if the \emph{$d$-vertex extremal length} of $\Gamma(v)$, \label{extlenght}
$$d\op{-VEL}(\Gamma)(v):= \sup_{m}\inf_{\gamma \in \Gamma(v)} \frac{\op{Length}_{m}(\gamma)^d}{\|m\|^d_{d}}$$
is infinite. It is easily seen that this definition does not depend upon the choice of $v\in V$. This natural extension of VEL parabolicity from \cite{HS} can be found earlier in \cite{BS09}.
\begin{rqe} In the context of bounded degree graphs, $2$-parabolicity is equivalent to recurrence  of the simple random walk on the graph, see \cite{HS} and the references therein. In general, $2$-VEL  is closely related to discrete conformal structures such as circle packings and square tilings, see \cite{BSsquaretiling,BrSmStTu40,HS,Ssquaretiling}.
\end{rqe}

\subsection{Limit of graphs}
\label{limits}
A rooted graph $(G = (V,E), o \in V)$ is isomorphic to $(G' = (V',E'), o'\in V'
)$ if there is a graph-isomophism of $G$ onto $G'$ which takes $o$ to $o'
$. We can define (as introduced in \cite{BeSc01}) a distance $\Delta$ on the space of isomorphism classes of locally finite rooted graphs by setting
$$\Delta \big((G, o), (G', o')\big) =\Big(1+\sup \left \{k :\op{Ball}_{G}(o, k) \mbox{ isomorphic to } \op{Ball}_{G'}(o',k)\right \}\Big)^{-1},$$
where $\op{Ball}_{G}(o, k)$ is the closed combinatorial ball of radius $k$ around $o$ in $G$ for the graph distance. In this work, limits of graph should be understood as referring to $\Delta$. It is easy to see that the space of isomorphism classes of rooted graphs with maximal degree less than $M$ is compact
with respect to $\Delta$. In particular every sequence of random rooted graphs of degree bounded by $M$ admits
weak limits.

\begin{defi} \label{unbiaised} A random rooted graph $(G,o)$ is \emph{unbiased} if $(G,o)$ is almost surely finite and conditionally on $G$, the root $o$ is uniform over all vertices of $G$.\end{defi}

We are now ready to state our main result. The case $d=2$ is \cite[Theorem 1]{BeSc01}.

\begin{thm} \label{dpara} Let $M \geq 0$ and $d \in \{2,3, ...\}$. Let $(G_{n},o_{n})_{n\geq0}$ be a sequence of unbiased random rooted graphs such that, almost surely, for all $n\geq 0$, $G_{n}$ is $M$-uniform in dimension $d$. If $(G_{n},o_{n})$ converges in distribution towards $(G,o)$ then $G$ is almost surely $d$-parabolic.
\end{thm}
 Applications of Theorem \ref{dpara} will be discussed in Section 4. 

\section{Proof of Theorem \ref{dpara}}
We follow the structure of the proof of \cite[Theorem 1]{BeSc01}: \begin{enumerate}
\item We first construct a limiting random packing whose tangency graph \emph{contains} the limit of the finite graphs.
\item The main step consists in showing that this packing has at most one accumulation point (for the centers) in $\R^d$ , almost surely.
\item Finally we conclude by quoting  a theorem relating packing in $\R^d$ and $d$-parabolicity.
\end{enumerate}
 Let $(G_{n},o_{n})_{n\geq0}$ be a sequence of unbiased, $M$-uniform in dimension $d$, random rooted graphs converging to a random rooted graph $(G,o)$. Given $G_{n}$, let $P_{n}$ be a deterministic $M$-uniform packing of $G_{n}$ in $\R^d$. We can assume that $o_{n}$ is independent of $P_{n}$. \\

Suppose that $C\subset\R^d$ is a finite set of points (in the application below, $C$ will be the set of centers of balls in $P_{n}$).
When $w\in C$, we define its { isolation radius} as
$\rho_w:=\inf\bigl\{|v-w|:v\in C\setminus\{w\}\bigr\}$.
Given $\delta\in(0,1),\,s>0$ and $w\in C$, following \cite{BeSc01} we say that $w$
is \emph{$(\delta,s)$-{supported}} if
in the ball of radius $\delta^{-1}\rho_w$, there are more than
$s$ points of $C$ outside of every ball
of radius $\delta\rho_w$; that is, if
$$
\inf_{ p\in\R^d}\,\,
\Bigl|C\cap \op{Ball}_{\R^d}(w,\delta^{-1} \rho_w)\setminus \op{Ball}_{\R^d}(p,\delta\rho_w)
\Bigr|\ge s\,.
$$

\begin{center}
\includegraphics[height=7cm]{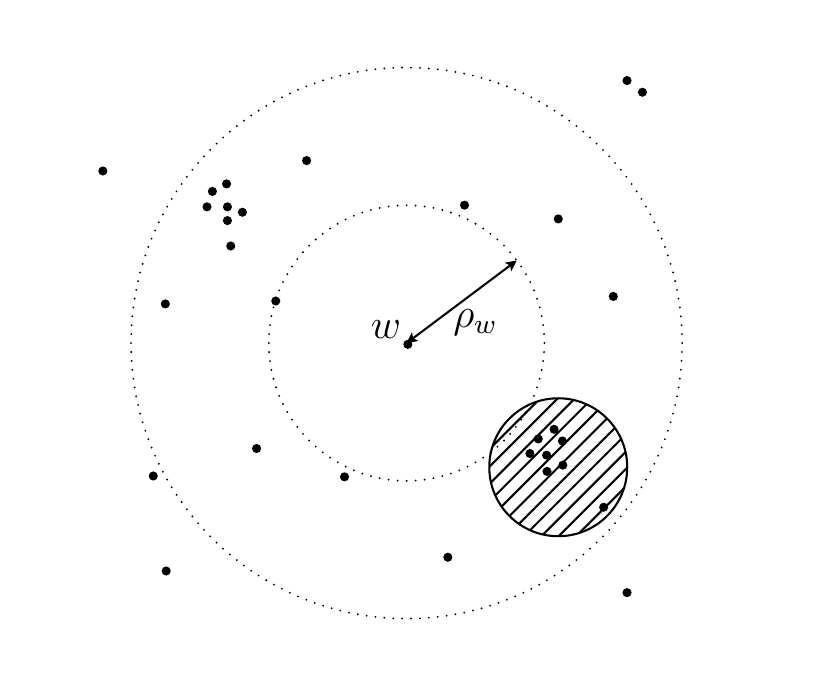}\\
\textsc{Fig.\,1:} Illustration of the definition of $(\delta,s)$-supported. Here, the point $w$ is $(0.5,7)$-supported
\end{center}

\begin{lem}[{\cite{BeSc01}}]\label{supported}
Let $d\geq 2$. For every $\delta\in(0,1)$ there is a constant
$c(\delta,d)$ such that
for every finite set $C\subset\R^d$ and every $s \ge 2$ the set
of $(\delta,s)$-supported points in $C$
has cardinality at most $c(\delta,d) |C|/s$.
\end{lem}
 Lemma 2.3 in \cite{BeSc01} deals with the case $d=2$, but the proof when $d\geq 2$ is the same  and is therefore omitted.  \\

Now, thanks to this lemma and to the fact that the point $o_{n}$ has been chosen independently of the packing $P_{n}$, for any $\delta>0$ and any $n \geq 0$, the probability that the center of the ball $B_{o_{n}}$ is $(\delta,s)$-supported in the centers of $P_{n}$ goes to $0$ as $s \to \infty$. Let $\tilde{P}_{n}$ be the image of $P_{n}$ under a linear mapping so that the ball $B_{o_{n}}$ is the unit ball in $\R^d$. Since the definition of $(\delta,s)$-supported is invariant under dilations and translations, we have \begin{eqnarray}
\mathbb{P}\big( 0 \mbox{ is }(\delta,s)\mbox{-supported in the centers of }\tilde{P}_{n}\big)  & \underset{s \to \infty}{ \longrightarrow} 0. \label{fix2}\end{eqnarray}
Let $\tilde{\mathbf{P}}_{n}$ be the union of the spheres of the packing $\tilde{P}_{n}$ and $\tilde{\mathbf{C}}_{n}$ be the union of the centers of the spheres of $\tilde{P}_{n}$. By definition, $\tilde{\mathbf{P}}_{n}$ and  $\tilde{\mathbf{C}}_{n}$ are random closed subsets of $\mathbb{R}^d$. The topology of Hausdorff convergence on every compact of $\R^d$ is a compact topology for closed subsets of $\R^d$. Hence, we can assume that along a subsequence we have the following convergence in distribution 
\begin{eqnarray}
\big((G_{n},o_{n}),\tilde{\mathbf{P}}_{n},\tilde{\mathbf{C}}_{n}\big) & \underset{n \to \infty}{\longrightarrow} \big((G,o),\mathbf{P},\mathbf{C}\big), \label{ascv}
\end{eqnarray}
related to $\Delta$ for the first component and to the Hausdorff convergence on every compact of $\R^d$ for the second and third ones. Without loss of generality we can suppose that there is no need to pass to a subsequence and by Skorhokhod representation theorem that the convergence \eqref{ascv} is almost sure.

\begin{pro} The random closed set $\mathbf{P}$ is almost surely the closure of a sphere packing in $\R^d$ whose centers have at most one accumulation point in $\mathbb{R}^d$. Furthermore, the tangency graph associated to $\mathbf{P}$ almost surely contains $(G,o)$ as a subgraph.
\end{pro}
\proof We begin with the second claim of the proposition.  By definition of $\tilde{P}_{n}$ we know $\mathbf{P}$ contains the unit sphere of $\mathbb{R}^d$ that corresponds to $o \in G$. Since the packings $\tilde{P}_{n}$ are $M$-uniform, any vertex neighbor of $o_{n}$ in $G_{n}$ corresponds to ball in the packing whose radius is in $[M^{-1},M]$ and tangent to the unit ball of $\mathbb{R}^d$. This property passes to the limit and by \eqref{ascv} we deduce that any neighbor of $o$ in $G$ corresponds to a sphere of $\mathbf{P}$ of radius in $[M^{-1},M]$ and tangent to the unit sphere of $\mathbb{R}^d$. A similar argument shows that $\mathbf{P}$ almost surely contains tangent spheres whose tangency graph contains $G$. Note that in the set $\mathbf{P}$ new connexions can occur (non tangent spheres in $\tilde{P}_{n}$ can become tangent at the limit).\\
The first part of the proposition reduces to showing that $\mathbf{C}$ almost surely has at most one accumulation point in $\mathbb{R}^d$. We argue by contradiction and we suppose that with probability bigger than $\varepsilon$, there exists two accumulation points $A_{1}$ and $A_{2}$ in $\mathbf{C}$ such that $|A_{1}-A_{2}| \geq \varepsilon$ and $|A_{1}|,|A_{2}| \leq \varepsilon^{-1}$. This implies, by \eqref{ascv}, that for any $s \geq 0$ with a probability asymptotically bigger than $\varepsilon$ the point $0$ is not $(\varepsilon/2, s)$-supported in $\tilde{\mathbf{C}}_{n}$. Which contradicts \eqref{fix2}.
\endproof

Since every subgraph of a $d$-parabolic graph is itself $d$-parabolic (obvious from the definition), the following extension of \cite[Theorem 3.1 (1)]{HS} together with the last proposition enables us to finish to proof of the Theorem \ref{dpara}.
\begin{thm}[{\cite[Theorem 7]{BS09}}] \label{BS2} Let $G$ be a graph of bounded degree. If $G$ is packable in $\R^d$ and if the packing has finitely many accumulation points in $\R^d$, then $G$ is $d$-parabolic.
\end{thm}
\begin{rqe} In order to be totally accurate, the $d$-parabolicity notion defined in \cite{BS09} corresponds to the definitions of Section \ref{dpot} when the function $m$ is defined on the edges of the graph. But these two notions easily coincide in the bounded degree case.
\end{rqe}

\section{Geometric applications} \label{geom}
\subsection{Isoperimetric inequalities and alternative}
If $W$ is a subset of a graph $G$, we recall that $\partial W$ is the set of vertices not in $W$ but neighbor with some vertex in $W$. We begin with an isoperimetric consequence of $d$-parabolicity which is an extension of \cite[Theorem 9.1(1)]{HS}. The proof is similar.
\begin{pro}\label{geo} Let $G=(V,E)$ be a locally finite, infinite, connected graph. Let $o \in V$, and $g : \R_{+} \to \R_{+}^*$ be some nondecreasing function. \begin{enumerate}[(1)]
\item  Suppose that $G$ is $d$-parabolic. If for every finite set $W$ containing $o \in W$, we have $|\partial W| \geq g(|W|)$ then
\begin{equation} \label{sum} \sum_{n=1}^\infty g(n)^{-\frac{d}{d-1}} = \infty. \end{equation}
\item If $g$ satisfies (\ref{sum}) and if $|\partial W_{k}| \leq g(|W_{k}|)$, for $(W_{k})_{k\geq0}$ defined recursively by
$$ W_{0}= \{o\} \mbox{ and } W_{k+1}= W_{k}\cup \partial W_{k} \mbox{ for }k \geq 0,$$
then $G$ is $d$-parabolic.
\end{enumerate}
\end{pro}
\proof  We know by assumption that $d\op{-VEL}(\Gamma(o)) = \infty.$ This implies that we can find functions $m_{i} : V \to \R_{+}$ such that $\|m_{i}\|_{d} = 2^{-i}$ and $\inf_{\gamma \in \Gamma(o)} \op{Length}_{m_{i}}(\gamma) \geq 1$. Hence $m := \sum_{i=0}^\infty m_{i}$ defines a function on $V$ such that $$\|m\|_{d} \leq 1 \mbox{ and } \inf_{\gamma \in \Gamma(o)} \op{Length}_{m}(\gamma)=\infty.$$
Without loss of generality we will suppose that $m(v) >0$ for all vertices $v\in V$.  The function $m \in \mathbb{L}^d(V)$ defines a distance on $V\times V$ by setting $$\op{d}_{m}(v,v'):=\inf \{ \op{Length}_{m}(\gamma), \gamma: v \to v' \}.$$ The idea is to explore the graph $G$ in a continuous manner according to $\op{d}_{m}$ and to use the isoperimetric inequality provided by $g$. For each $v \in V$ let $$I_{v}:= [\op{d}_{m}(o,v)-m(v),\op{d}_{m}(o,v)].$$ For $h \in \R_{+}$, we define $s_{v}(h):= \frac{\op{Leb}(I_{v} \cap [0,h])}{m(v)}$. Intuitively, water flows in the graph $G$ starting from $o$, $m(v)$ is the time that water needs to wet $v$ before flowing to its neighbors. A vertex $v\in V$ begin to get wet at $h= \min I_{v}$ and is completely wet at $h=\max I_{v}$. The function $s_{v}(h)$ represents the percentage of water in $v$. We set $s(h):= \sum_{v\in V}s_{v}(h).$ Since $\op{d}_{m}(o,\infty)=\infty$, for every $h \in \R_{+}$ there are only finitely many $v\in V$ such that $s_{v}(h) \ne 0$ and then $s(h)$ is piecewise linear. We denote $W_{h}:=\{v\in V, h \geq \max I_{v}\}$ the set of vertices that are totaly wet at time $h$ and $G_{h}:= \{v\in V, h \in I_{v}\}$ the set of vertices that are getting wet at time $h$. Clearly $G_{h} = \partial W_{h}.$ Let $$ f(x) = \min \left( g\left( \frac{x}{2}\right),\frac{x}{2} \right).$$
If $|G_{h}| \geq s(h)/2$ then \begin{equation} \label{iso}|G_{h}| \geq f(s(h)),\end{equation} otherwise $|G_{h}| < s(h)/2,$ then the number of completely wet vertices is at least $s(h)/2$ and consequently $|G_{h}| \geq g(s(h)/2)$.  Thus (\ref{iso}) always holds.\\
At points where $h\mapsto s(h)$ is differentiable we have $$\frac{ds}{dh}(h) = \sum_{v \in G_{h}}s'_{v}(h)= \sum_{v\in G_{h}}\frac{1}{m(v)}.$$
Writing $1= m(v)^{(d-1)/d}m(v)^{-(d-1)/d}$ and using H\"older inequality with $p=d$ we get
$$ \left( \sum_{v\in G_{h}} 1\right) \leq \left(\sum_{v \in G_{h}} \frac{1}{m(v)} \right)^{\frac{d-1}{d}} \left( \sum_{v \in G_{h}} m(v)^{d-1}\right)^{1/d},$$ and thus using (\ref{iso}):
$$\frac{ds}{dh}(h) \geq \frac{|G_{h}|^{\frac{d}{d-1}}}{\left(\sum_{v\in G_{h}}m(v)^{d-1}\right)^{\frac{1}{d-1}}} \geq \frac{f(s(h))^{\frac{d}{d-1}}}{\left(\sum_{v\in G_{h}}m(v)^{d-1}\right)^{\frac{1}{d-1}}},$$
$$\mbox{ therefore }\frac{ds}{f(s(h))^{\frac{d}{d-1}}}\geq \frac{dh}{\left(\sum_{v\in G_{h}}m(v)^{d-1}\right)^{\frac{1}{d-1}}}.$$
Integrating for $0<a<h<b<\infty$ and using H\"older with $p=d$ we get
$$\int_{s(a)}^{s(b)} \frac{ds}{f(s)^{\frac{d}{d-1}}}ds \geq \int_{a}^b \frac{dh}{\left(\sum_{v\in G_{h}}m(v)^{d-1}\right)^{\frac{1}{d-1}}} \geq \frac{(b-a)^{d/(d-1)}}{\left(\int_{a}^b \left(\sum_{v\in G_{h}}m(v)^{d-1}\right) dh\right)^{1/(d-1)}}.$$
Remark that $\int_{0}^\infty \left(\sum_{v\in G_{h}}m(v)^{d-1}\right) dh = \sum _{v\in V}m(v)^d < \infty$, and that $s(b) \to \infty$ when $b \to \infty$. We conclude that the integral of $f(.)^{-\frac{d}{d-1}}$ diverges and the same conclusion holds for $g(.)^{-\frac{d}{d-1}}$. Since $g(.)$ is non-decreasing, a comparison series-integrale ends the proof of the first part of the proposition.\\
For the second part, set $n_{k}= |W_{k}|$ and define for $N \in \N^*$ a function $m: V \to \R_{+}$ on $G$ by
$$ m(v) = \left\{ \begin{array}{ll} g(n_{k})^{-\frac{1}{d-1}} & \mbox{for }v\in \partial W_{k} \mbox{ and } k \leq N, \\ 0 &\mbox{otherwise.}
\end{array} \right.
$$ Then we have  $\inf \{ \op{Length}(\gamma) : \gamma \in \Gamma(o)\}\geq \sum_{k=0}^N g(n_{k})^{-\frac{1}{d-1}}$ and
$$ \|m\|_{d}^{d} \leq \sum_{k=0}^N \frac{|\partial W_{k}|}{g(n_{k})^{d/(d-1)}} \leq \sum _{k=0}^N g(n_{k})^{-\frac{1}{d-1}}.$$
By definition of the extremal length, it suffices to show that
$\sum_{k=0}^\infty g(n_{k})^{-\frac{1}{d-1}}=\infty.$ Note that $n_{k+1} \leq n_{k} + g(n_{k})$, thus by monotonicity of $g$, we obtain
$$ \frac{1}{g(n_{k})^{\frac{1}{d-1}}} \geq \frac{1}{n_{k+1}-n_{k}} \sum_{n=n_{k}}^{n_{k+1}-1} \frac{1}{g(n)^{\frac{1}{d-1}}} \geq \sum_{n=n_{k}}^{n_{k+1}-1} \frac{1}{g(n_{k})}\frac{1}{g(n)^{\frac{1}{d-1}}} \geq \sum_{n=n_{k}}^{n_{k+1}-1} \frac{1}{g(n)^{d/(d-1)}}.$$
Which implies $\displaystyle \sum_{k=0}^\infty g(n_{k})^{-\frac{1}{d-1}} \geq \sum_{n_{0}}^\infty g(n)^{-d/(d-1)} = \infty.$  \endproof

Let us recall the definition of the Cheeger constant of a infinite graph $G$:
$$\op{Cheeger}(G) := \inf \left\{ \frac{|\partial W|}{|W|}: W \subset G, |W| < \infty \right\}.$$
The following corollary generalizes a theorem regarding planar graphs indicated by Gromov and proved by several authors. See  Bowditch \cite{Bow95}  for a very short proof and references for previous proofs.
\begin{cor} Let $G$ be an infinite locally finite connected graph which admits a $M$-uniform packing in $\R^d$. Then we have the following alternative:
\begin{itemize}
\item either $G$ has a positive Cheeger constant,
\item or they are arbitrarily large subsets $W$ of $G$ such that
$$ |\partial W| \leq |W|^{\frac{d-1}{d}+o(1)}.$$
\end{itemize}
\end{cor}
\proof Let $G$ be a infinite connected graph which is the tangency graph of a $M$-uniform packing in $\R^d$ (in particular $G$ has  bounded degree). If $\op{Cheeger}(C)=0$, then we can find a sequence of subsets $A_{i} \subset G$ such that $$\frac{|\partial A_{i}|}{|A_{i}|}\underset{i\to \infty}{\longrightarrow} 0.$$ We associate to $A_{i}$ the random unbiased graph $(A_{i},o_{i})$ where $o_{i}$ is uniform over the vertices of $A_{i}$. By compactness (see the discussion before Definition \ref{unbiaised}),  along a subsequence we have the weak convergence for $\Delta$
$$ (A_{i},o_{i}) \underset{i \to \infty}{\overset{(d)}{\longrightarrow}} (A,o),$$
where $(A,o)$ is almost surely $d$-parabolic. We assume that there is no need to pass to a subsequence. Therefore the sequence of rooted random graphs $(A_{i},o_{i})_{i \geq 1}$ satisfies all the hypotheses of Theorem \ref{dpara}, in particular $(A,o)$ is almost surely $d$-parabolic. Let $\varepsilon>0$. By Proposition \ref{geo}, almost surely, there exists a (random) subset $W\subset A$ containing $o \in W$ and satisfying $$|\partial W| \leq |W|^{\frac{d-1}{d}+\varepsilon}.$$ We claim that the set $W$ and its boundary are contained in $G$. Indeed for any $k \geq 0$, the bounded degree assumption combined with the fact that $\frac{|\partial A_{i}|}{|A_{i}|} \to 0$ imply that
$$ \mathbb{P}\left( o_{i} \mbox{ is at a graph distance less than $k$ from }\partial A_{i}\right) \underset{i \to \infty}{\longrightarrow} 0.$$ Hence, almost surely for any $k\geq0$, the ball of radius $k$ around $o$ in $A$ is a subgraph of some $A_{i}$'s and thus of $G$. This finishes the proof of the corollary.
\endproof

\subsection{Non existence of $M$-uniform packing}
As a consequence of the last corollary, the graph $\Z^{d+1}$ cannot be $M$-uniform packed in $\R^d$ for some $M\geq 0.$ This is a weaker result compared to \cite{BS09} where it is shown that $\Z^{d+1}$ cannot be sphere packed in $\R^d$ using non-existence of bounded non constant $d$-harmonic functions on $\Z^d$. 

The \emph{parabolic index} of a graph $G$ (see \cite{Soardi-Yamasaki}) is the infimum of all $d\geq0$ such that $G$ is $d$-parabolic (with the convention that $\inf \varnothing = \infty$). For example, Maeda \cite{Maeda} proved that the parabolic index of $\Z^d$ is $d$. It is easy to see that the parabolic index of a regular tree is infinite, leading to the following consequence.
\begin{cor} Let $G_{n}$ be a deterministic sequence of finite graphs. If there exists $f(n) \underset{n \to \infty}{\longrightarrow} \infty$ and $k \in \{2,3,...\}$ such that
$$ \frac{\# \{ v \in G_{n}, \op{Ball}_{G_{n}}(v,f(n)) = \ k\mbox{-regular tree up to level }f(n) \}}{|G_{n}|} \underset{n \to \infty}{\longrightarrow} 1,$$
then for all $M \geq 0$, $G_{n}$ eventually cannot be $M$-uniform packed in $\R^d$.
\end{cor}
\proof Note that any unbiased weak limit of $G_{n}$ is the $k$-regular tree and apply Theorem \ref{dpara}.
\endproof
That is, if for a sequence of $k$-regular graphs, $ k >2$, the girth grows  to infinity  then only finitely many of the graphs
can be $M$-uniform  packed in any fixed dimension. The same holds if the limit is some other nonamenable graph.

\section{Open problems}

Several necessary conditions are provided in this paper for a graph to be ($M$-uniform) packed in $\R^d$. The first two questions are related to existence of packable graphs in $\R^d$.
\begin{ques} \begin{enumerate}   \item Find necessary and sufficient conditions for a graph to be ($M$-uniform) packable in $\R^d$.
\item Exhibit a natural family of graphs which are ($M$-uniform) packable in $\R^d$.
\item Show that the number of tetrahedrangulations  in $\R^3$ with $n$ vertices grows to infinity.
\end{enumerate}
\end{ques}

\begin{ques}
It is of interest to understand  what is the analogue of packing of a graph and the results above in the context of Riemannian manifolds. Does packable in the discrete context of graphs is analogous to conformally flat?
\end{ques}

\begin{ques} Show that the Cayley graph of Heisenberg group $\mathbf{H}_{3}(\Z)$ generated by
$$ A = \left(\begin{array}{ccc} 1&1&0\\ 0&1&0\\0&0&1 \end{array}\right) \mbox{ and } B = \left(\begin{array}{ccc} 1&0&0\\ 0&1&1\\0&0&1 \end{array}\right),$$
is not packable in $\R^d$ though is known to be  $4$-parabolic, see e.g. \cite{S94}.
\end{ques}

 The two following questions deal with the geometry of the accumulation points (of centers) of packing in $\R^d$. 

\begin{ques} Does there exist a graph $G$ packable in $\R^d$ in two manners $P_{1}$ and $P_{2}$ such that the set of accumulation points in $\mathbb{R}^d \cup\{\infty\}$ for $P_{1}$ is a point but not for $P_{2}$ ?
\end{ques}

\begin{ques}[\cite{BS09}] Show that any packing of $\Z^3$ in $\R^3$ has at most one accumulation point in $\mathbb{R}^d \cup \{\infty\}$.
\end{ques}


\begin{ques}[Parabolicity for edges] What is left of Theorem \ref{dpara} in the context of edge parabolicity (where the function $m$ of Section \ref{dpot} is defined on the edges of the graph) without the bounded degree assumption ? For instance, is it the case that every  limit of unbiased random planar graphs is 2-edge-parabolic (which means SRW is recurrent) ?
\end{ques}

\begin{ques}[Sub-diffusivity] Let $G$ be a $d$-parabolic graph. Consider $(S_{i})_{i\geq0}$ a simple random walk on $G$. Do we have
$$
\liminf_{n\to \infty} \frac{\op{d}_{\op{gr}}(S_{0},S_{n})}{\sqrt{n}} \ne \infty \quad?$$
\end{ques}

\begin{ques}[Mixing time] Let $G$ be a finite graph packable in $\R^d$ with bounded degree. Show that mixing time is bigger than $C_d \op{diameter}(G)^2.$ In particular the planar $d=2$ case is still open.
\end{ques}

\noindent\textbf{Acknowledgments:} Part of this work was done during visit of the second author to Weizmann Institute. The
 second author thanks his hosts for this visit. We are indebted to Steffen Rohde and to James T. Gill for pointing out several inaccuracies in a first version of this work. Thanks for the referee for a very useful report.

\bibliographystyle{abbrv}

\end{document}